\documentclass[12pt, titlepage]{article}
\input{xypic}
\usepackage{amsfonts}
\usepackage{amssymb}
\usepackage{latexsym}
\def\qed{{\hbadness=10000\hfill\ \vbox{\hrule height.09ex
   \hbox{\vrule width.09ex height1.55ex depth.2ex \kern1.8ex
   \vrule width.09ex height1.55ex depth.2ex}\hrule height.09ex}\break
   \bigskip}}

\newtheorem{theorem}[equation]{Theorem}
\newtheorem{fact}[equation]{Fact}

\newtheorem{proposition}[equation]{Proposition}
\newtheorem{example}[equation]{Example}
\newtheorem{lemma}[equation]{Lemma}
\newtheorem{corollary}[equation]{Corollary}

\newenvironment{remark}{\refstepcounter{equation}\noindent
{\bf Remark \thesection.\arabic{equation}} }{$\blacksquare$}
\makeatletter
\@addtoreset{equation}{section}
\makeatother
\newcommand{\proof}{\noindent {\bf Proof~}}

\renewcommand{\pmod}[1]{\;(\text{mod } #1)}
\newcommand{\PGL}{{\mathrm{PGL}}}	
\renewcommand{\char}{\text{char}\,}	
\newcommand{\Z}{{\mathbb Z}} 

\title{Ovals and Hyperovals in Desarguesian Nets}

{\normalsize \author{by David A. Drake and Kevin Keating\\
Department of Mathematics \\
University of Florida \\
Gainesville, FL 32611 \\
USA \\[.2cm]
\date{}
$\quad$\\
{\tt dad@math.ufl.edu}\\
{\tt keating@math.ufl.edu}}
\date{}
\begin{document}

\maketitle

\begin{abstract}
We determine the Desarguesian planes which hold $r$-nets with ovals
and those which hold $r$-nets with hyperovals for every $r \le 7$.
\end{abstract}

\section{Introduction} \label{S:intro}

A {\it net} is an incidence structure $\Sigma$ whose blocks are
partitioned into three or more {\it parallel classes} so that each
class consists of two or more blocks which partition the points of
$\Sigma$ and so that any two blocks not in the same parallel
class have a unique point of intersection.  The blocks of a net
$\Sigma$ are referred to as {\em lines}.
One says that $\Sigma$ is an {\it $r$-net} (or a net of {\it
degree} $r$) if the number of parallel classes is $r < \infty$,
a net of {\it order} $n$ if some line has cardinality
$n < \infty$.  If $\Sigma$ has order $n$, each line and
parallel class has cardinality $n$. 

If $\Sigma$ is a net obtained from a net $\Pi$ by deleting zero or
more complete parallel classes of $\Pi$, one says that $\Sigma$ is
{\it embedded} in $\Pi$ or that $\Sigma$ {\it extends} to $\Pi$ or
that $\Pi$ {\it holds} $\Sigma$.  We often take $\Pi$ to be an {\it
affine plane}; i.\,e., a net in which each pair of points is joined by
a line.  If $\Pi^*$ is the projective plane obtained by adjoining a
single line (denoted by $\ell_{\infty}$) to an affine plane $\Pi$, we
use the same language for $\Pi^*$ as for $\Pi$.  We call $\Sigma$ a
{\it Desarguesian} net if it is held by a Desarguesian plane.

A set $S$ of points in a net $\Sigma$ is called an {\it arc} if each 
pair of points of $S$ is joined by a line of $\Sigma$ that intersects
$S$ only in those two points.  A {\it $k$-arc} is an arc of
cardinality $k$.  We define an {\it oval} and a {\it hyperoval} in an
$r$-net to be an $r$-arc and an $(r + 1)$-arc, respectively.  In this
paper, we investigate the problem of determining the pairs $(r,\Pi)$
for which $\Pi$ is a Desarguesian plane that holds an $r$-net with
oval or hyperoval.  We solve both problems for $r \le 7$.

Ovals in $r$-nets held by {\it finite} Desarguesian projective planes
have been investigated by G. J. Simmons \cite{Simm}.  Simmons
describes such an oval as an $r$-set ``sharply focused on''
$\ell_{\infty}$.  He utilized these ovals to construct geometry-based
secret sharing schemes.  We thank the referee who brought this
interesting paper to our attention. 

Let $\Pi$ be the Desarguesian affine plane coordinatized by a division
ring $D$.  We prove that $\Pi$ holds a 5-net with oval if and only if
$D$ contains a root of $x^2 + x - 1$ (Theorem~\ref{t8.2} below, a
result obtained by Hirschfeld for finite $D$); we prove 
that $\Pi$ holds a 7-net with oval if either $D$
contains $GF(2^k)$ for some $k \ge 3$ or $D$ contains a root of
$x^3 - x^2 - 2x + 1$ (Theorem~\ref{t5.3}).  We prove that $\Pi$ holds
a 6-net with oval if and only if either char $D \ne 2, 3$ or $D$
properly contains $GF(4)$ (Theorem~\ref{t4.6}).

Let $H$ be a hyperoval of an $r$-net $\Sigma$. Then every point of $H$
is joined to the remaining points of $H$ by a line of each parallel
class of $\Sigma$, so there are no tangents to $H$.  Thus, a parallel
class of $\Sigma$ partitions the points of $H$ into $(r + 1)/2$
subsets of size 2, so $r$ must be odd.  Let $\Pi$ be a Desarguesian
affine plane coordinatized by a division ring $D$.  We prove that
$\Pi$ holds a 5-net with hyperoval if and only if $D$ contains $GF(4)$
(Theorem~\ref{t8.3}) and that $\Pi$ holds a $7$-net with hyperoval if 
and only if char $D = 2$ and $\vert D \vert \ge 8$
(Theorem~\ref{t5.1}). 

Let $S$ and $T$ both be sets of points or both be ordered sets of
points of an affine plane $\Pi$.  We say that $S$ and $T$ are {\it
affinely equivalent} or {\it $\Pi$-equivalent} if there is a
collineation $\phi$ of $\Pi$ with $(S)\phi = T$.  For a Desarguesian
affine plane $\Pi$, we prove that all ovals (hyperovals) of $r$-nets
held by $\Pi$ are $\Pi$-equivalent for each $r \le 6$
(Proposition~\ref{t3.3} and Theorems~\ref{t8.2}, \ref{t8.3},
\ref{t4.6}).  Hyperovals in 7-nets are not all affinely equivalent,
however (see, Remark~\ref{t5.7}).

In Section 6, we present a number of non-existence results
(Corollaries~\ref{t6.2}, \ref{t6.4}, \ref{t6.6}).  In
Corollaries~\ref{t3.13} and \ref{t6.10}, we determine the values of
$r$ for which there exist Desarguesian or general $r$-nets of order
$n$ with ovals or hyperovals for small $n$.

It is an open question whether there exist any Desarguesian nets of
characteristic not equal to 2 with hyperovals, or any Desarguesian nets
of degree not equal to $2^k\pm1$ with hyperovals.

\section{Constructions} \label{S:construct}

\begin{fact} \label{t2.1}
(i) Removing any point of a hyperoval of a net $\Sigma$ produces an
oval of $\Sigma$. 
\par
(ii) If a net $\Pi$ holds a net $\Sigma$ with hyperoval, adjoining an
additional parallel class of $\Pi$ to $\Sigma$ produces a net with
oval. 
\end{fact}

If $D$ is a division ring, we shall write $\Pi(D)$ and $\Pi^*(D)$
for the Desarguesian affine and projective planes coordinatized by
$D$.  The point set of $\Pi(D)$ is $D \times D$; the lines of $\Pi(D)$
are the sets of points $(x,y)$ which satisfy an equation of one of the
forms $x = c$, $y = mx + b$ with $b, c, m$ in $D$.  The points of
$\Pi^*(D)$ are the homogenous triples $\langle x,y,z \rangle :=
\{(xt,yt,zt):t\in D \setminus \{0\}\}$ with $(x,y,z)\not=(0,0,0)$.   
One embeds $\Pi(D)$ in $\Pi^*(D)$ by identifying the point
$(x,y)\in\Pi(D)$ with the point $\langle x,y,1 \rangle\in \Pi^*(D)$
and adjoining an {\it ideal line}
$\ell_{\infty}$ consisting of the {\it ideal points} $\langle x,y,0
\rangle$.  In particular, the point $\langle 1,m,0 \rangle$ is added
to the lines of $\Pi(D)$ of slope $m$, and the point $\langle 0,1,0
\rangle$ is added to the lines of infinite slope.  An oval or
hyperoval of a net $\Sigma$ held by $\Pi(D)$ is called a {\it
subgroup} or {\it coset oval} or {\it hyperoval} if it is a subgroup
or a coset of a subgroup of the additive group $D\times D$.

Multiplication in a division ring $D$ need not be commutative, so we
adopt the convention that $b/c$ and $\frac{b}{c}$ both denote
$bc^{-1}$.  We frequently utilize the following observation.

\begin{fact} \label{t2.11}
Let $b$ be an element of a division ring $D$.  Then the sub-division
ring of $D$ generated by $b$ is commutative.
\end{fact}

\begin{fact} \label{t2.2}
If $\Pi(D)$ holds an $r$-net $\Sigma$ with coset hyperoval $H$,
then $\vert D \vert > r$ or ${(\vert D \vert , r) = (2,3)}$.
\end{fact}

\proof Let $\vert D \vert = q < \infty$.  If $r = q + 1$, then $\vert
H \vert = q + 2$;  so Lagrange's Theorem implies that $(q + 2) \,
\vert \, q^2$, hence that $q = 2$ and $r = 3$.  If $r = q$, one
obtains the contradiction $(q + 1) \, \vert \, q^2$.~$\blacksquare$

\begin{proposition} \label{t2.4}
Let $D$ be a division ring of characteristic 2, and let $k$ be a
positive integer such that $2^k \le \vert D \vert$.  Then $\Pi(D)$
holds a $(2^k - 1)$-net with a subgroup hyperoval $H$.  One may choose
$H$ so that each pair of parallel secants of $H$
intersect $H$ in an affine subplane of order 2.
\end{proposition}

\proof If $D$ is finite, $D$ is a field;  if $D$ is
infinite, $D$ contains an
infinite field \cite[(13.10)]{L}.  In both cases, $D$ contains a
field $F$ with ${|F| \ge 2^k}$.  Let $S$ be a subgroup of order
$2^k$ of the additive group of $F$, and set ${G= \{ (x,x^2) \, \vert
\, x \in S \}}$.  Since $F$ is a field of characteristic 2, $G$ is a
subgroup of $F \times F$, and $\vert G \vert = \vert S \vert = 2^k$.
Each point $(x,x^2)$ in $G$ is joined to the remaining $2^k-1$ points
$(x+y,x^2+y^2)$ of $G$ by lines whose slopes are precisely the $2^k-1$
elements of $S \setminus \{ 0 \}$.  Thus, no three points of $G$ lie
in a common line, and $G$ is a hyperoval in the net which consists of
the parallel classes of $\Pi(F)$ with slopes in $S \setminus \{ 0 \}$.
(Note that if $S = F$ and $F$ is finite, the union of $G$ and the
appropriate two ideal points forms a ``regular'' hyperoval or
``complete conic'' of $\Pi^*(F)$, cf.\ \cite[p.\,32]{C} or
\cite[p.\,79]{E}.) 

Suppose that the secants $(a,a^2)(b,b^2)$ and $(c,c^2)(d,d^2)$ have a
common slope $a+b = c+d$.  Then $a + b + c + d = 0$, and the six
secants determined by the indicated four points lie in the three
parallel classes with slopes $a+b$, $a+c$, $a+d$.~$\blacksquare$

\begin{corollary} \label{t2.3}
Suppose that $D$ is a division ring of characteristic 2 with
$\vert D \vert \ge 2^k$.  Then $\Pi(D)$ holds a $2^k$-net with
subgroup oval and a $(2^k-1)$-net with oval. 
\end{corollary}

\proof Apply Proposition~\ref{t2.4} and Fact~\ref{t2.1}. 
$\blacksquare$

\begin{proposition} \label{t2.5}
Suppose that the division ring $D$ contains a field $F$
isomorphic to $GF(q)$.  Then

$\,\,\,\,(i) \,\,$ $\Pi(D)$ holds an $r$-net with oval for $r = q -
1$, $q$, $q + 1$; 

$\,\, (ii) \,\,$ if $q$ is even, $\Pi(D)$ holds a $(q + 1)$-net with
hyperoval; 

$(iii) \,\,$ if $q$ is even and $F \ne D$, $\Pi(D)$ holds a $(q +
2)$-net with oval. 
\end{proposition}

\proof It is well known (see \cite[pp.\,567--569]{BJL} or
\cite{H}) that the Desarguesian projective plane $\Pi^*(F)$ 
of order $q$ contains an oval $O$ and that $\Pi^*(F)$ contains a
hyperoval $H$ if $q$ is even.  The plane $\Pi(D)$ contains $\Pi(F)$
which is obtained from $\Pi^*(F)$ by removing a line
$\ell_{\infty}$.  If $r$ denotes $\vert \ell_{\infty} \setminus O
\vert$, the $r$ parallel classes of $\Pi(D)$ determined by the points
of $\ell_{\infty} \setminus O$ form an $r$-net with oval $O \setminus
\ell_{\infty}$.  Taking $\ell_{\infty}$ to be a secant, tangent or
exterior line to $O$ in $\Pi^*(F)$, one obtains $r$ equal to $q-1$,
$q$, $q+1$, respectively.  If $q$ is even and $\ell_{\infty}$ is
disjoint from $H$, the lines of $\Pi(D)$ with slopes in $F \cup \{
\infty \}$ are a $(q+1)$-net with hyperoval $H$.  The truth of
(iii) follows from (ii) and Fact~\ref{t2.1}\,(ii).~$\blacksquare$

\begin{proposition} \label{t2.7} Let $F$ be a field contained
in a division ring $D$; let $r\geq3$; and let $\zeta$ be a primitive 
$r$-th root of unity in the algebraic closure of $F$.
Assume that $\zeta + \zeta^{-1} \in F$.  Then $\Pi(D)$
holds an $r$-net with oval.
\end{proposition}

\proof For $k \ge 2$ we have
\[\zeta^k + \zeta^{-k} = (\zeta + \zeta^{-1})^k +
\sum_{i=1}^{k-1}\,c_i (\zeta^i + \zeta^{-i})\]
with $c_i\in F$,
so $\zeta^k + \zeta^{-k} \in F$  for all $k$.  Let
\[P_k=\left(\zeta^k + \zeta^{-k},
\sum_{i = -k+1}^{k-1} \zeta^i\right)\in\Pi(D),\]
and let  $O=\{P_k \, \vert \, 0 < k \le r\}$.  For $0 < k, \, \ell \le
r$, $k \ne \ell$, the line $P_kP_{\ell}$ has slope $(\zeta^{k + \ell}
+ \zeta)/[(\zeta^{k + \ell} - 1)(\zeta - 1)]$.  In particular,
$P_kP_{r-k}$ has slope $\infty$ if $k \ne r, \, r/2$.  On the other
hand, if $P_kP_{\ell}$ has finite slope $m$ then
\[\zeta^{k+\ell} = \frac{(\zeta-1)m + \zeta}{(\zeta-1)m - 1}.\]
Thus, each $P_k$ is joined to the other
points of $O$ by lines with $r - 1$ distinct slopes; so
$O$ is a set of cardinality $r$, and
no three points of $O$ are collinear.
Since the lines joining points of $O$ have only $r$ distinct
slopes, $O$ is an oval in an
$r$-net held by $\Pi(D)$.~$\blacksquare$

\begin{proposition} \label{t2.9}
Let $r\geq3$, let $p$ be an odd prime not dividing $r$, and
let $\zeta$ be a primitive $r$-th root of unity in
the algebraic closure of $GF(p)$.  Then the smallest extension
of $GF(p)$ containing $\zeta+\zeta^{-1}$ is $GF(p^b)$, where
$b$ is the order of $p$ in $(\Z/r\Z)^{\times}/(\pm1)$.
\end{proposition}

\proof The field $GF(p)(\zeta)$ is an extension of
$GF(p)(\zeta+\zeta^{-1})$ whose degree $d$ is 1 or 2.
Let $a$ be the order of $p$ in $(\Z/r\Z)^{\times}$; then
$a=[GF(p)(\zeta):GF(p)]$, and hence
$[GF(p)(\zeta+\zeta^{-1}) : GF(p)] = a/d$.
We have $d=2$
if and only if $a$ is even and $\zeta+\zeta^{-1}$ is fixed
by the map $x\longrightarrow x^{p^{a/2}}$.  Thus, $d = 2$ precisely when 
$\zeta^{p^{a/2}}+\zeta^{-p^{a/2}}=\zeta+\zeta^{-1}$ or,
equivalently, when $(\zeta^{p^{a/2}}-\zeta^{-1})(1-\zeta^{1-p^{a/2}})=0$;
i.\,e., when $p^{a/2}\equiv-1\pmod r$.  Thus $a/d$ is the order
of $p$ in $(\Z/r\Z)^{\times}/(\pm1)$.~$\blacksquare$

\medskip

The truth of the following proposition is asserted by Simmons in
\cite[p.\,204]{Simm}.  Simmons illustrates his construction method for
values of $q$ up to 19.

\begin{proposition} \label{t2.10} Suppose that the division ring $D$
contains a field isomorphic to $GF(q)$ with $q$ odd and $q \equiv \pm
1 \pmod r$, $r \ge 3$.  Then $\Pi(D)$ holds an $r$-net with oval.
\end{proposition}

\proof Apply Propositions~\ref{t2.7} and~\ref{t2.9}.~$\blacksquare$

\section{Ovals and Hyperovals in Desarguesian Nets of Small Degree}
\label{S:small} 

Every collineation of a Desarguesian plane $\Pi(D)$ can be expressed
as the composition of a {\it linear collineation} of the form 

\[
(x , y) \,\, \longrightarrow \,\, (ax + by , cx + dy) 
\,\, + \,\, (e , f) \quad \text{for some } \,\, 
\left( \begin{array}{cc} a & b \\ c & d \end{array} \right)\in
GL_2(D)
\]
and some $e$, $f \in D$ with a collineation of the form
$(x,y) \longrightarrow ( \, x^{\alpha}, y^{\alpha} \, )$ for some
$\alpha \in$ Aut$(D)$.

A set $Q$ of four points of a net $\Sigma$ is called a {\it quad} if
$Q$ is a quadrangle (no three points of $Q$ are collinear) and if two
(or more) parallel classes of $\Sigma$ each contain two secants to
$Q$.  An ordered set $(A, B, C, E)$ of points of $\Sigma$ is said to
be an {\it ordered quad} if $AB$ and $CE$ are parallel lines of
$\Sigma$ and $AE$ and $BC$ are parallel lines of $\Sigma$.

\begin{proposition} \label{t3.1}
(well known) (i) Let $m_0$, $m_1$, $m_{\infty}$ be distinct elements
of $D \cup \{ \infty \}$ and let $A,B$ be points of $\Pi(D)$ such
that $AB$ has slope $m_1$.  Then there is a linear collineation 
mapping $A$ to $(0,0)$, $B$ to $(1,1)$, and taking lines of slopes
$m_0$, $m_1$, $m_{\infty}$, respectively, into lines of slopes 0, 1,
$\infty$.  

(ii) Every ordered quad of $\Pi(D)$ is affinely equivalent to
$\{(0,0)$, $(1,0)$, $(1,1)$, $(0,1)\}$.  Hence, a quad of $\Pi(D)$ is
an affine subplane of $\Pi(D)$ if and only if char $D = 2$.
\end{proposition}

\proof It is well known that $\PGL_3(D)$, acting by left matrix
multiplication on the points of $\Pi^*(D)$, represented as column
vectors, is transitive on the ordered quadrangles of $\Pi^*(D)$.  In 
the arguments below, $\langle 1, \infty, 0 \rangle$ should be
interpreted as $\langle 0, 1, 0 \rangle$.

To prove (i), let $\phi^*$ be an element of $\PGL_3(D)$ which
maps the quadrangle
\[(A, B,\langle 1, m_0, 0 \rangle,\langle 1, m_{\infty},
0\rangle )\]
to the quadrangle
\[(\langle 0,0,1 \rangle, \langle 1,1,1
\rangle, \langle 1,0,0 \rangle, \langle 0,1,0 \rangle ).\]  The
collineation $\phi^*$ fixes the line $\ell_{\infty}$ and thus may be
represented as left multiplication by a non-singular $(3 \times
3)$-matrix whose third row is $(0 \, 0 \, 1)$.  Thus $\phi^*$
induces a linear collineation $\phi$ of $\Pi(D)$.  Clearly, $\phi$
maps $A$ to $(0,0)$ and $B$ to $(1,1)$.  Since $\phi^*$ maps  
$\langle 1, m_1, 0 \rangle$ to $\langle 1, 1, 0 \rangle$, $\phi$ has
the desired action on parallel classes of $\Pi(D)$.

To prove (ii), let $\mu^*$ be an element of $\PGL_3(D)$ which
maps a quadrangle ${(A, B, C, E)}$ to the quadrangle
\[( \langle 0,0,1
\rangle, \langle 1,0,1 \rangle, \langle 1,1,1 \rangle, \langle 0,1,1
\rangle ).\]
The collineation $\mu^*$ maps $AB \cap CE$ to
$\langle 1,0,0 \rangle$ and $AE \cap BC$ to
$\langle 0,1,0 \rangle$ and, hence,
fixes the line $\ell_{\infty}$. Thus $\mu^*$ induces a linear
collineation $\mu$ of $\Pi(D)$ with the specified action on
the points $A$, $B$, $C$, $E$.~$\blacksquare$

\medskip

\begin{corollary} \label{t3.2}
If \,$\Pi(D)$ holds an $r$-net $\Sigma$ with coset hyperoval
$H$, then $\Pi(D)$ holds an $r$-net $\Sigma'$ with subgroup hyperoval
$H'$.  Furthermore, for each ${0\leq i\leq 3}$ such that
$i \le \vert D \vert - r + 1$, we may choose $\Sigma'$
to exclude any $i$ of the parallel classes of lines of slope 0, 1,
$\infty$ and to include the other $3 - i$ classes. 
\end{corollary}

\proof Let $A$ be any point of $H$.  Choose $m_0$, $m_1$,
$m_{\infty} \in D \cup \{ \infty \}$ so that lines of slope $m_h$ are
in $\Sigma$ if and only if one desires lines of slope $h$ to be in
$\Sigma'$.  Apply Proposition~\ref{t3.1}\,(i).~$\blacksquare$

\begin{fact} \label{t3.5}
Let $O$ be an oval in an $r$-net $\Sigma$.  If $r$ is even, at least
$r/2$ parallel classes of \,$\Sigma$ contain (exactly) $r/2$
secants to $O$.  If $r$ is odd, each parallel class of $\Sigma$
contains $(r - 1)/2$ secants to $O$.
\end{fact}

\proof There are $r(r - 1)/2$ secants to $O$ and at most $r/2$
secants in each of the $r$ parallel classes of $\Sigma$.~$\blacksquare$

\begin{proposition} \label{t3.3}
Let $D$ be a division ring.  Then $\Pi(D)$ holds a 3-net with oval;
$\Pi(D)$ holds a 4-net with oval if and only if $\vert \, D \, \vert
\ne 2$; and  $\Pi(D)$ holds a 3-net with hyperoval if and only if $D$
has characteristic 2.  All ovals of $r$-nets held by $\Pi(D)$ are
affinely equivalent for $r = 3$ and for $r = 4$, and all hyperovals of
3-nets held by $\Pi(D)$ are affinely equivalent.
\end{proposition}

\proof Any set of three non-collinear points in $\Pi:=\Pi(D)$ is an
oval in a 3-net held by $\Pi$.  Assume that $\Pi$ holds a 3-net with a
hyperoval $H$.  Disjoint secants to $H$ are parallel lines of $\Pi$;
so $H$ is an affine subplane of $\Pi$.  By Proposition~\ref{t3.1}\,(ii), one
may assume that $H = \{(0,0) , (0,1) , (1,0) , (1,1) \}$.  Since $H$
is an affine plane, $\{0,1\}\subset D$ is a field; so $D$ has
characteristic 2, and $H$ is a subgroup of $D \times D$.  Conversely,
let
$S = \{(0,0) , (0,1) , (1,0) , (1,1) \}$.
If $\char D \ne 2$, $S$ is an oval in the 4-net whose lines have slopes 0,
1, $-1$, $\infty$.  If $\char D = 2$, $S$ is a hyperoval in the 3-net
$\Sigma$ whose lines have slopes 0, 1, $\infty$.  In this case, one
obtains a 4-net with oval by Fact~\ref{t2.1}\,(ii) unless $\vert \, D
\, \vert = 2$.

Clearly, every hyperoval of a 3-net is a quad; by Fact~\ref{t3.5},
every oval of a 4-net is a quad;  it is easily seen that every oval of
a 3-net is contained in a quad.  Hence the last three claims
follow from Proposition~\ref{t3.1}\,(ii).~$\blacksquare$

\begin{proposition} \label{t3.50} (Hirschfeld \cite[Lemma 7.1.2]{H})  
Let $D$ be a finite field.  If $\Pi(D)$ holds a 5-net with
oval, then $D$ contains a root of $x^2 - x - 1$.
\end{proposition}

An oval $O$ of a 5-net of $\Pi(D)$ is said to be in {\it standard
position} if $O = \{(1,1), (1,0), (0,0), (0,b), (c,b) \} =: O_{b,c}$
for some $b$, $c$ in $D \setminus \{0, 1\}$.  The following
proposition removes the finiteness assumption of
Proposition~\ref{t3.50}. 

\begin{proposition} \label{t3.6}
Let $D$ be a division ring.  Then every oval $O$ of every 5-net
$\Sigma$ of $\Pi(D)$ is affinely equivalent to an oval in standard
position.  If $O_{b,c}$ is an oval in standard position, $c = b + 1$;
and $b$ is a root of $x^2 + x - 1$.  (We denote $O_{b,c}$ by $O_b$.) 
\end{proposition}

\proof By Fact~\ref{t3.5}, each of the five parallel classes of
$\Sigma$ contains 
two secants to $O$.  If $O$ were the union of a quad $S$ and a point
$P$, one of the five parallel classes would contain none of the six
secants to $S$, so would contain at most one secant to $O$.  This
contradiction of Fact~\ref{t3.5} proves that $O$ contains no quad.  From
Proposition~\ref{t3.1}\,(i), one sees that $O$ is affinely
equivalent to some $O_{b,c}$.

The point $(c,b)$ is on no secant of slope $\infty$, so must be on a
secant of each of the other four slopes.  The secant of slope 1
containing $(c,b)$ must be $(1,0)(c,b)$, so $c = b + 1$.  Then
$(0,b)(1,1) \parallel (0,0)(b+1,b)$.  Thus, $1 - b = b/(b+1)$; so $0 =
b^2 + b - 1$.~$\blacksquare$

\begin{theorem} \label{t8.1} (Gordon - Motzkin \cite[(16.4)]{L})
Let $a_1, \dots , a_n$ be elements of a division ring $D$.
Then each root of the polynomial $f(x) = (x - a_1) \dots (x - a_n)$
is conjugate to some $a_i$.
\end{theorem}

The second assertion of the following theorem generalizes and sharpens
an assertion of Simmons \cite[p.\,197]{Simm} that the ovals of 5-nets
of $\Pi(GF(11))$ are projectively equivalent.

\begin{theorem} \label{t8.2}
The Desarguesian plane $\Pi(D)$ holds a 5-net with oval if and only if
$D$ contains a root $b$ of $x^2 + x - 1$.  All ovals of 5-nets held by
$\Pi(D)$ are affinely equivalent.
\end{theorem}

\proof If $\Pi(D)$ holds a 5-net with oval, Proposition~\ref{t3.6}
asserts that $D$ contains a root $b$ of $x^2 + x - 1$.  Conversely,
assume that $D$ contains an element $b$ with $b^2 = 1 - b$.  Then $b
\ne 0, \pm 1$; so $b + 1 \ne 0$, 1.  The set $O_b$ has two secants
of each of the slopes 0, 1, $\infty$, $1 - b$ and $-b$; so $O_b$ is an
oval of a 5-net held by $\Pi(D)$.  

By Proposition~\ref{t3.6}, each oval $O$ is affinely equivalent to an
oval $O_d$.  Let $O_b$ be a fixed oval in standard position, and let
$\phi$ be the collineation given by right multiplication by the matrix
\[
\left( \begin{array}{cc} -b & -1-b \\ 1+b & 1+b \end{array} \right) \,
.
\]
Then $\phi$ maps $O_b$ to $O_c$ where $c = -b - 1$; so $O_b$ and $O_c$
are affinely equivalent.  By Theorem~\ref{t8.1}, there is an inner
automorphism $\alpha$ of $D$ with $d^{\alpha} = b$ or $d^{\alpha} =
c$.  Then the 
collineation $(x,y) \longrightarrow (x^{\alpha},y^{\alpha})$ maps
$O_d$ to $O_b$ or $O_c$.~$\blacksquare$

\begin{corollary} \label{t3.7}
(i) If $\char D \not= 2$, $\Pi(D)$ holds a 5-net with oval if and
only if $D$ contains an element $a$ such that $a^2=5$.

(ii) If $\char D = p > 0$, then $\Pi(D)$ holds a 5-net with
oval if and only if $p=5$ or $p \equiv \pm 1 \pmod{5}$ or $D$ contains
a field isomorphic to $GF(p^2)$.   
\end{corollary}

\proof  Part (i) follows from the quadratic formula.  To prove
(ii), first note that $x^2 + x - 1 = 0 $ has solutions in $GF(4)$ but 
no solution in
$GF(2)$.  Furthermore, by quadratic reciprocity (see, for
instance, \cite[p.\,7, Theorem 6]{S}), if $p>2$ then $x^2=5$
has a solution in $GF(p)$ if and only if $p=5$ or
$p\equiv\pm1\pmod5$.  Part (ii) now follows using
Fact~\ref{t2.11}.~$\blacksquare$

\begin{theorem} \label{t8.3} The plane $\Pi(D)$ holds a 5-net with
hyperoval if and only if $D$ is a division ring which contains a field
isomorphic to GF(4).  Let $b \in D$ be a root of $x^2 + x - 1$.
All hyperovals of all 5-nets held by $\Pi(D)$ are affinely equivalent
to the point set 
$H_b := \{ (1,1), (1,0), (0,0), (0,b), (b+1,b), (b+1,1) \}$.
\end{theorem}

\proof The sufficiency of the condition that $D$ contain $GF(4)$ is a
special case of Proposition~\ref{t2.5}\,(ii).  Conversely, suppose that
$\Pi(D)$ holds a 5-net $\Sigma$ with hyperoval $H$.  By
Fact~\ref{t2.1}\,(i) and Proposition~\ref{t3.6}, $H = \{ (1,1)$, $(1,0)$,
$(0,0)$, $(0,b)$, $(b+1,b)$, $(c,d) \}$ for some $b$, $c$, $d$ in
$D$ with $b^2 + b - 1 = 0$.  The secants $(c,d)(1,1)$ and
$(c,d)(b+1,b)$ must have respective slopes $0$ and $\infty$, so $(c,d)
= (b+1,1)$.  The only secant through $(b+1,1)$ which can have the
slope $-b$ of secant $(1,0)(0,b)$ is $(b+1,1)(0,0)$, so $-b =
1/(b+1)$.  Then $-b$ is also the slope $(b-1)/b$ of $(b+1,b)(1,1)$.
Thus, $b^2 + b + 1 = 0 = b^2 + b - 1$.  It follows that char $D = 2$
and that $\{ 0, 1, b, b+1 \} \subseteq D$ is a field.

Suppose that $H_b$ and $H_c$ are hyperovals in 5-nets held by
$\Pi(D)$.  By Theorem~\ref{t8.2}, there is a collineation $\phi$ of
$\Pi(D)$ with $(O_b)\phi = O_c$; so $(H_b)\phi = H_c$.~$\blacksquare$ 

\begin{proposition} \label{t3.9}
If $D$ contains a root of $x^3 - x^2 - 2x + 1$ then
the Desarguesian plane $\Pi(D)$ holds a 7-net with oval.
\end{proposition}

\proof Suppose that $b^3 - b^2 - 2b + 1 = 0$ with $b\in D$.
Then $b \ne 0, \pm 1$, so the set \[O = \{ (0,0), (1,0), (0,-1),
(1,b-1), (1-b, -1), (b^2-b, b-1), (1-b, b-b^2) \}\]
consists of seven distinct elements.
Let $\mathfrak{S}$ denote the multiset consisting of the 21
secants $PQ$ with $P, Q \in O$.  The secant
$(1,0)(b^2-b, b-1)$ has slope $(b-1)/(b^2-b-1) = b$,
and the secant $(1-b, -1)(b^2-b, b-1)$ has slope $b/(b^2-1) = b - 1$.
It is now easy to see that $\mathfrak{S}$ contains three secants
of each of the following slopes: 1, 0, $b$, $b-1$.  One observes that
the map $(x,y) \longrightarrow (-y,-x)$ fixes $O$ and maps lines of
slope $m$ to lines of slope $1/m$ for each $m$.  Thus,
$\mathfrak{S}$
also contains three secants of each of the slopes $\infty$, $1/b$,
$1/(b-1)$.  One checks that the seven listed slopes are distinct.
Thus, $\mathfrak{S}$ has three lines of each of the slopes 1, 0, $b$,
$b-1$, $\infty$, $1/b$, and $1/(b-1)$,
so no secant meets $O$ in more
than two points.  Thus, $O$ is an oval in a 7-net held by
$\Pi(D)$.~$\blacksquare$

\begin{theorem} \label{t5.3} The plane $\Pi(D)$ holds a 7-net with
oval if and only if either $D$ contains a field isomorphic
to $GF(2^k)$ with $k \ge 3$ or $p(x) := x^3 - x^2 -  2x + 1$
has a root in $D$.  
\end{theorem}

\proof If $D$ contains $GF(2^k)$ with $k \ge 3$ one obtains a
7-net with oval in $\Pi(D)$ by Corollary~\ref{t2.3}.  If
$p(x)$ has a root in $D$, one obtains a 7-net with oval in
$\Pi(D)$ by Proposition~\ref{t3.9}.  The converse is proved
in \cite{D2}.~$\blacksquare$

\begin{lemma} \label{t3.10}
The polynomial $x^3 - x^2 - 2x + 1$ has a root in $GF(p)$ if
and only if $p=7$ or $p \equiv \pm 1 \pmod{7}$.  
\end{lemma}

\proof If $p=7$ then $x=-2$ is a root of $x^3-x^2-2x+1$.  If
$p\not=7$ let $\zeta$ be a primitive 7th root of 1 (in an
algebraic closure of $GF(p)$).  Then $-\zeta-\zeta^{-1}$,
$-\zeta^2-\zeta^{-2}$, $-\zeta^3-\zeta^{-3}$ are the
roots of $x^3-x^2-2x+1$.  We have
$-\zeta^{pi}-\zeta^{-pi}=-\zeta^i-\zeta^{-i}$ if and only if
$\zeta^{-pi}(\zeta^{(p-1)i}-1)(\zeta^{(p+1)i}-1)=0$.  Thus if
$p\equiv\pm1\pmod7$ then all the roots of $x^3-x^2-2x+1$ lie
in $GF(p)$, while if $p\not\equiv\pm1\pmod7$, none of them
do.~$\blacksquare$

\begin{corollary} \label{t5.4} 
If $\char D = p \ne 0$, then $\Pi(D)$ holds a 7-net with oval if and
only if one of the following conditions holds:

$\,\,\,\, (i) \,\, D$ contains a field isomorphic to $GF(2^k)$
for some $k \ge 3$; 

$\,\, (ii) \,\, p = 7$ or $p\equiv\pm 1\pmod{7}$;

$ (iii) \,\, D$ contains a field isomorphic to $GF(p^3)$.
\end{corollary}

\proof Apply Theorem~\ref{t5.3} and Lemma~\ref{t3.10}.~$\blacksquare$

\section{Desarguesian 6-Nets with Ovals}\label{S;six}

For the following four lemmas, we take $S$ to be a 6-set, a {\it
(parallel) class} on $S$ to be a collection of 2-subsets of $S$ which 
partitions $S$, $\Omega_r$ to be a set of $r$ mutually disjoint
parallel classes $\Pi_1, \dots, \Pi_r$ on $S$.  We say that $(S,
\Omega_r)$ and $(S', \Omega'_r)$ are {\it isomorphic} if there is a
bijection $\phi: S \longrightarrow S'$ and a permutation $\sigma$ of
$\{ 1, \dots, r \}$ with $(\Pi_i)\phi = \Pi'_{(i)\phi}$ for $1 \le i
\le r$ where $\Omega'_r = \{ \Pi'_1, \dots, \Pi'_r \}$.

\begin{lemma} \label{t4.1}  (P. Cameron \cite[Theorem 4.7 (ii)]{CA}) 
Each of $(S, \Omega_2)$ and $(S, \Omega_5)$ is unique up to
isomorphism.  The automorphism group of $(S, \Omega_5)$ is
5-transitive on $\Omega_5$.
\end{lemma}

We take $S$ to be the set of six cells $(i,j)$ of a $3 \times 3$
matrix for which $i + j \ne 4$.  Since $(S, \Omega_2)$ is unique, we
may take the rows and columns of this matrix to be the lines of
$\Pi_1$ and $\Pi_2$, respectively, for each $(S, \Omega_r)$ with $r
\ge 2$.  The following lemma is easy to verify.

\begin{lemma} \label{t4.2}
In $(S, \Omega_3)$, the lines of $\Pi_3$ are represented by the
entries of one of the following matrices.
$$
\underbrace
{\pmatrix
{1 & 2 & * \cr
3 & * & 1 \cr
* & 3 & 2}}_{A_1},
\qquad \,\,
\underbrace
{\pmatrix
{1 & 3 & * \cr
2 & * & 3 \cr
* & 1 & 2}}_{A_2},
\qquad \,\,
\underbrace
{\pmatrix
{1 & 2 & * \cr
2 & * & 3 \cr
* & 3 & 1}}_{A_3},
\qquad \,\,
\underbrace
{\pmatrix
{1 & 2 & * \cr
3 & * & 2 \cr
* & 3 & 1}}_{A_4}.
$$
\end{lemma}

A structure $(S, \{ \Pi_1, \dots, \Pi_r \})$ is said to be {\it
maximal} if there is no class $\Pi_{r+1}$ which is disjoint from
$\Pi_1 \cup \dots \cup \Pi_r$.

\begin{lemma} \label{t4.3}\,(i) In $(S, \Omega_5)$, the lines of
$\Pi_3$, $\Pi_4$, $\Pi_5$ may be represented by the entries of $A_1$,
$A_2$, $A_3$. 

(ii) In $(S, \Omega_3)$, the class $\Pi_3$ may be represented by the 
entries of $A_3$ or $A_4$;  if by $A_4$, then $(S, \Omega_3)$ is
maximal.  

(iii)  The automorphism group of $(S, \Omega_3)$ is 3-transitive on
$\Omega_3$. 
\end{lemma}

\proof The classes of lines induced by $A_1$, $A_2$, $A_3$ are
mutually disjoint, so conclusion (i) follows from Lemma~\ref{t4.1}.
The first assertion of conclusion (ii) and the $A_3$-case of
conclusion (iii) follow from conclusion (i) and Lemma~\ref{t4.1}.

Suppose that $\Pi_3$ is represented by $A_4$.  By Lemma~\ref{t4.2},
the only contenders for a class $\Pi_4$ disjoint from $\Pi_1 \cup
\Pi_2 \cup \Pi_3$ is one represented by $A_i$ for some $i \le 3$.  
None of these is disjoint from $\Pi_3$, so the second assertion of
conclusion (ii) is valid.  The map $\phi: S \longrightarrow S, (i,j)
\longrightarrow (j,i)$ induces the permutation $(\Pi_1,\Pi_2)$ on
$\Omega_3$.  If $\mu$ is the permutation of $S$ given by
$((1,2),(2,1),(3,3))$, then $\mu$ induces the permutation $(\Pi_1,
\Pi_2, \Pi_3)$ on $\Omega_3$.  Thus, conclusion (iii) is also valid in
the $A_4$-case.~$\blacksquare$  

\begin{lemma} \label{t4.4}
Let $(S, \Omega_3)$ be contained in the plane $\Pi(D)$.  Then there 
are $a, b \in D \setminus \{ 0, 1 \}$ and a linear collineation of
$\Pi(D)$ mapping the elements of $S$ to 
$$
\left( \begin{array}{ccc}
(0,0) & (a,0) & *\\
(0,b) & * & (1,b)\\
\ast & (a,1) & (1,1)
\end{array} \right).
$$
\end{lemma}

\proof By Proposition~\ref{t3.1}\,(ii), we may assume that the points 
represented
by cells $(1,1)$, $(1,3)$, $(3,1)$ and $(3,3)$ are the points with
coordinates $(0,0)$, $(1,0)$, $(0,1)$ and $(1,1)$,
respectively.~$\blacksquare$

\begin{theorem} \label{t4.6}
The Desarguesian affine plane $\Pi(D)$ holds a 6-net $\Sigma$ with
oval $O$ if and only if either $\char D \ne 2$, $3$ or $D$ properly
contains a field isomorphic to $GF(4)$.  If $\char D \ne 2$, $\Sigma$ 
induces only three complete parallel classes on $O$; if $\char D = 2$,
$O$ is a hyperoval of a 5-net held by $\Sigma$.  In both cases, all
ovals of all 6-nets held by $\Pi(D)$ are affinely equivalent.
\end{theorem}

\proof The point set $O = \{(0,0), (1,0), (0,1), (2,1), (1,2), (2,2)
\}$ has three secants of each of the slopes 0, 1, $\infty$ and two 
secants of each of the slopes $1/2$, 2, $-1$.  If $\char D \ne 2,3$,
these six values are distinct, so $\Pi(D)$ holds a 6-net with oval
$O$.  On the other hand, by Theorem~\ref{t8.3} and
Fact~\ref{t2.1}\,(ii), if $D$ properly contains a field isomorphic to
$GF(4)$, then $\Pi(D)$ contains a 6-net with oval. 

Assume, conversely, that the Desarguesian affine plane $\Pi(D)$ holds
a 6-net $\Sigma$ with oval $O$.  By Fact~\ref{t3.5}, $\Sigma$ contains
(at least) three parallel classes which contain three secants to $O$.
By Lemma~\ref{t4.3}\,(ii), the secants of these three classes may be
represented by matrix $A_3$ or matrix $A_4$.  Furthermore, we may
assume that the coordinates of the points of $O$ are given by the
matrix of Lemma~\ref{t4.4}. 

We first treat the case of $A_4$.  By Lemma~\ref{t4.3}\,(iii),
$\Sigma$ induces exactly three complete parallel classes on $O$ and,
hence, three partial classes consisting of two secants each.  Note
that the lines of one of the complete parallel  
classes represented by $A_4$ have slope $1=b/(1-a)=(1-b)/a$,
so $a+b=1$.  Therefore the slopes of the lines joining the
points $(a,0)$, $(0,b)$, $(1,1)$ are $a$, $1/b$, $-b/a$, and
the slopes of the lines joining the points $(0,0)$, $(1,b)$,
$(a,1)$ are the reciprocals of these numbers.
Since these six lines lie in three parallel classes, and $a$,
$1/b$, $-b/a$ are distinct, we must have
$\{a,1/b,-b/a\}=\{1/a,b,-a/b\}$.
If $a=1/a$, then since $a\not=1$ we get $a=-1$,
$\char D\not=2$ and $b=2$.  If $a=b$, then $2a=1$; so $\char D\not=2$
and $a=b=1/2$.  If $a=-a/b$, then $b=-1$ and $a=2\ne0$, so
$\char D\not=2$.  In all three cases
we get $\{a,1/b,-b/a\}=\{-1,2,1/2\}$, and hence $\char D\not=3$. 

Let $O_1$, $O_2$, $O_3$ be the respective ovals produced by taking $a
= 1/a$, $a = b$, and $a = -a/b$.  Then $O_1$ is mapped to $O_2$ and
$O_3$ under the collineations induced by the respective matrices 
\[
\left( \begin{array}{cc}
0 & 1\\
1 & 0\\
\end{array} \right)
\qquad
\text{and}
\qquad
\left( \begin{array}{cr}
0 & -1/2\\
1 & 3/2\\
\end{array} \right),
\]
so the three ovals are affinely equivalent.

Now suppose that the secants of three parallel classes of $\Sigma$ are
represented by 
$A_3$.  Then $(a,0)(0,b)$ and $(1,b)(a,1)$ have slope 1, so $b = -a$
and ${1 - b = a - 1}$.  Thus, $\char D = 2$ and $b = a$.  Let $M$ denote
the multiset $\{(1+b)/b,1/(1+b),b\}$, and let $M'$ denote the multiset
of reciprocals $\{b/(1+b),1+b,1/b\}$.  There are three secants each of
slopes $0$, $1$, $\infty$;  and the slopes of the remaining six
secants are the multiset $S := M \cup M'$.  If $S$ consists of two
distinct numbers, each of multiplicity 3, then $\Sigma$ contains a
5-net with hyperoval.  Thus, Theorem~\ref{t8.3} implies that $D$
contains $GF(4)$ and that all ovals of all 6-nets held by $\Pi(D)$ are
affinely equivalent.

On the other hand, suppose that $S$ contains more than two numbers.
If any element of $M$ matches any element of $M'$, one obtains a
contradiction.  Thus, each of $M$ and $M'$ must contain
at least two distinct numbers; and we obtain the contradiction that $O$ has
secants of at least seven distinct slopes.~$\blacksquare$

\section{Desarguesian 7-Nets with Hyperovals}
\label{V:seven}

\begin{theorem} \label{t5.1}
The plane $\Pi(D)$ holds a 7-net with hyperoval $H$ if and only if $D$
is a division ring of characteristic 2 with $\vert D \vert \ge 8$.  If
$\Pi(D)$ holds a 7-net with hyperoval, it holds a 7-net with
subgroup hyperoval which is the disjoint union of two affine
subplanes. 
\end{theorem}

\proof By Proposition~\ref{t2.4}, $\Pi(D)$ holds a 7-net with a
subgroup hyperoval which is the disjoint union of two affine subplanes
if $D$ is a division ring of characteristic 2 with $\vert D
\vert \ge 8$. 

Conversely, assume that $\Pi := \Pi(D)$ for some division ring $D$;
and suppose that $\Pi$ holds a 7-net $\Sigma$ with hyperoval $H$.  We
treat first the case that some subset of four points of $H$ is a
quad.  By Proposition~\ref{t3.1}\,(ii), we may assume that the eight points of 
$H$
receive the coordinates displayed in the matrix
$$
M \,\, = \,\,
\pmatrix
{(0,0) & (1,0) & * & * \cr
(0,1) & (1,1) & * & * \cr
* & * & (a,c) & (b,c) \cr
* & * & (a,d) & (b,d)}
$$
for some elements $a, b, c, d \in D \setminus \{ 0, 1 \}$ with $a \ne 
b$, $c$ and $d \ne b$, $c$.  Assume, by way of contradiction, that
$\char D \ne 2$. 
Then the slope of $(0,1)(1,0)$ is $-1 \ne 1$.  By permuting the last two
rows and the last two columns of the matrix $M$ as well as the symbols
$a$, $b$, $c$, $d$, we may assume that the three remaining secants of
slope 1 are $(0,1)(a,d)$, $(1,0)(b,c)$ and $(a,c)(b,d)$.  Then
$d - 1 = a$, $c = b - 1$ and $a - b = c - d$.  Adding the three
equations gives $2d = 2b$, and hence $d = b$.  This is a
contradiction, so $\char D = 2$ if the hyperoval contains a quad.

The remaining case is that $H$ contains no quads.  Here, we may
coordinatize $\Pi$ by $D$ as in the proof of Lemma~\ref{t4.4} so that
the eight points of $H$ receive the coordinates displayed in the
matrix   
$$
\pmatrix 
{(0,0) & (a,0) & * & * \cr
(0,c) & * & (b,c) & * \cr
* & (a,d) & * & (1,d) \cr
* & * & (b,1) & (1,1)}
$$
for some elements $a, b, c, d \in D \setminus \{ 0, 1 \}$ with $a \ne
b$ and $c \ne d$.  Let $\Pi_1$, $\Pi_2$, $\Pi_3$  be the parallel
classes which contain, respectively, the horizontal lines, the
vertical lines, and the lines of slope 1.  Since $\Pi_3$ does
not generate a quad in $H$ with either $\Pi_1$ or $\Pi_2$,
$\Pi_3$ cannot contain either of the secants $(a,0)(b,1)$,
$(0,c)(1,d)$.  Thus, the secants from $\Pi_3$ may be represented by
one of the following four matrices.
$$
\pmatrix
{1 & 2 & * & * \cr
2 & * & 4 & * \cr
* & 3 & * & 4 \cr
* & * & 3 & 1},
\pmatrix
{1 & 2 & * & * \cr
2 & * & 3 & * \cr
* & 3 & * & 4 \cr
* & * & 4 & 1},
\pmatrix
{1 & 2 & * & * \cr
3 & * & 2 & * \cr
* & 3 & * & 4 \cr
* & * & 4 & 1},
\pmatrix
{1 & 2 & * & * \cr
3 & * & 4 & * \cr
* & 4 & * & 2 \cr
* & * & 3 & 1}.
$$
If the lines of $\Pi_3$ are
given by the leftmost matrix, one sees that
\[1 = \frac{-c}{a} = \frac{b -a}{1 - d} = \frac{d - c}{1 - b}.\]
It follows that $2a = 0$ and, hence, that $2 = 0$.  Similar
computations yield the same conclusion for the other three
matrices, so we have $\char D=2$.  Since there are 7 slopes
we have $|D|+1\ge7$, and hence $|D|\geq8$.~$\blacksquare$

\begin{proposition} \label{t5.2} Let $H$ be a hyperoval of a 7-net held by
$\Pi(D)$.  If $H$ has no quads, then there is a field $F$
contained in $D$ and a linear collineation $\phi$ of $\Pi(D)$ such
that $(H)\phi \subseteq F \times F$.
\end{proposition}

\proof We may choose $\phi$ so that $(H)\phi$ is represented
as in the proof of Theorem~\ref{t5.1}.  Suppose,
for example, that the lines of $\Pi_3$ are given by the leftmost 
matrix displayed in the proof of Theorem~\ref{t5.1}.  Then $c
= a$ and $d = 1 + a + b$.  The secant $(0,0)(b,1)$ has slope $1/b
\ne 0, 1, \infty$.  There must be a secant of the same slope through
the point $(1,1)$, so
\[\frac{1}{b} \in \left \{\frac{1}{1+a},\,1+a,\,
\frac{1+a}{1+b},\,\frac{a+b}{a+1} \right\}.\]
If, for example, $1/b = (1+a)/(1+b)$,
then $b = 1/a$; so $b$ is an element of the field $F$
generated over $GF(2)$ by $a$.  Similarly $b$, hence also $c$ and $d$, are in
$F$ if $1/b$ is one of the other three slopes.  The argument is
similar for the other three matrices.~$\blacksquare$

\medskip

We observe in Example~\ref{t5.5} below that the conclusion of
Proposition~\ref{t5.2} does not hold for arbitrary hyperovals of
Desarguesian 7-nets and, in Example~\ref{t5.6}, that there are
hyperovals which satisfy the hypothesis of Proposition~\ref{t5.2}.

\begin{example} \label{t5.5} Let $D$ be a division ring of
characteristic 2 which is not commutative.  Then $\Pi(D)$ holds a
7-net with hyperoval $H$ such that $(H)\phi$ is not contained in $F
\times F$ for any collineation $\phi$ of $\Pi(D)$ and any field $F$
contained in $D$.   
\end{example}

\proof  Let $a$, $b$ be any two non-commuting elements of $D$.  Then
$a$, $b$, 1 are linearly independent over the prime field
contained in $D$.
Consider the point sets $G = \{(0,0), (0,1), (1,0), (1,1) \}$ and $H
= G \cup (G + (a,b))$.  Each point of $H$ is joined to the remaining
seven points of $H$ by secants of the seven different slopes
0, 1, $\infty$, $b/a$, $(b+1)/a$, $b/(a+1)$ and $(b+1)/(a+1)$.
Thus, no three points of $H$ are collinear, and
$H$ is a hyperoval of a 7-net held by $\Pi(D)$.

Suppose that there is a collineation $\phi$ of $\Pi(D)$ and a
field $F$ contained in $D$ such that $(H)\phi\subseteq F \times F$.
By Theorem~\ref{t3.1}\,(ii) there is a linear collineation $\psi$
of $\Pi(F)=F\times F$ such that $\phi\circ\psi$ fixes the
points $(0,0)$, $(0,1)$, $(1,0)$, $(1,1)$.  Then
$\phi\circ\psi$ is induced by an automorphism $\alpha$
of $D$, so we have $H=G\cup (G + (a^{\alpha},b^{\alpha}))$.
Since $a^{\alpha}$ doesn't commute with $b^{\alpha}$, this is
a contradiction.~$\blacksquare$

\begin{example} \label{t5.6} Let $F = GF(8)$.  Then
$\Pi(F)$ holds a 7-net with a hyperoval $H$ which has no quads.
\end{example}

\proof Let $b \in F$ be a root of $x^3+x+1$, and let
\[H = \{ (0,0), (1,1), (0,1/b), (1/b,0), (1,b), (b,1), (b,1/b), (1/b,b) \}.\]
The secants
through $(0,0)$ have the seven distinct slopes 1, $\infty$, 0, $b$,
$1/b = b^2+1$, $1/b^2 = b^2+b+1$ and $b^2$; as do the secants
through each of the other points of $H$.  Thus, $H$ is a hyperoval of
the 7-net with the seven listed slopes.  Since $\char F = 2$, one can
verify that $H$ contains no quads by verifying that $H$ contains no
affine subplanes of order 2.  If $S$ is an affine subplane of
$\Pi(F)$ contained in $H$, so is $H \setminus S$ (with the same three
parallel classes).  Thus, it suffices to verify that there is no
affine subplane $S$ of $\Pi(F)$ with $(0,0) \in S \subseteq H$.
Any affine subplane of $\Pi(F)$ containing $(0,0)$ is a
subgroup of $F\times F$ of order 4, and it is easily verified
that $H$ contains no such subgroup.
Therefore $H$ contains no quads.~$\blacksquare$

\bigskip

\begin{remark} \label{t5.7}
Let $\Sigma$ be a $(q - 1)$-net with hyperoval $H$ which
is held by the affine plane $\Pi(GF(q))$ obtained from a
projective plane $\Pi^*$ by deleting a line $\ell_{\infty}$.  Clearly
the union of $H$ with the obvious two points of $\ell_{\infty}$ is a
hyperoval $H^*$ of $\Pi^*$.  Removal of any secant $\ell'_{\infty}$
from $\Pi^*$ produces an affine plane $\Pi'$ isomorphic to $\Pi$ which
holds a $(q-1)$-net $\Sigma'$ with hyperoval $H' := H^* \setminus (H^*
\cap \ell'_{\infty})$.  When $q = 8$, all hyperovals of $\Pi^*$ are
projectively equivalent; i.\,e., are images of each other under
collineations of $\Pi^*$ (see \cite[Theorems 5.2.4, 9.2.3]{H}).  The
constructions of Examples~\ref{t5.5} and \ref{t5.6} demonstrate,
however, that they are not all affinely equivalent; i.e., the
collineation group of $\Pi^* = \Pi^*(GF(8))$ is not transitive on
pairs $(H^*, \ell_{\infty})$ where $H^*$ is a hyperoval of $\Pi^*$ and
$\ell_{\infty}$ is a secant of $H^*$.
\end{remark}

\section{Non-Existence Results}\label{V;non}

Let $O^n$ denote the set of all integers $r$ for which there exists an
$r$-net of order $n$ with an oval.  If $n$ is a prime power, let
$O^n_d$ denote the set of all integers $r$ for which there exists a
Desarguesian $r$-net of order $n$ with an oval.  If ``oval'' is
replaced by ``hyperoval'' in the preceding two definitions, we obtain
sets which we label $H^n$ and $H^n_d$.

\begin{theorem} \label{t6.1}
(B. Segre, see \cite[10.3.3, Cor.\ 2 to 10.3.3, Cor.\ 2 to 10.4.4]{H}).
Let $A$ be a $k$-arc in the Desarguesian projective plane $\Pi^*$ over
$GF(q)$.  Suppose that $q$ is
even and $k > q - \sqrt{q} + 1$ or that $q$ is odd and $k > q -
\sqrt{q}/4 + 7/4$.  Then the only points $X$ for which $A \cup \{ X
\}$ is a $(k + 1)$-arc of $\Pi^*$ are 

(i) ($q$ even) the $q + 2 - k$ points of $H \setminus A$ where
$H$ is the unique hyperoval containing $A$;

(ii) ($q$ odd) the $q + 1 - k$ points of $O \setminus A$ where
$O$ is the unique oval containing $A$.
\end{theorem}

\begin{corollary} \label{t6.2}
Define $f(q)$ to be $q - \sqrt{q}$ for even $q$ and $q - \sqrt{q}/4 +
3/4$ for odd $q$.

(i) If $q - 2 \ge r > f(q)$, then $r \notin O^q_d$.

(ii) Let $r > f(q) - 2$.  Then $r \notin H^q_d$ if $r \le q-3$, and
$r \notin H^q_d$ if $r \le q-2$ and $q$ is even.
\end{corollary}

\proof
Assume, by way of contradiction, that $\Pi := \Pi(GF(q))$ holds an
{$r\text{-net}$} $\Sigma$ 
with oval $O$ for some $r > f(q)$.  Adjoin a line $\ell_{\infty}$ to $\Pi$
to form the Desarguesian projective plane $\Pi^*$.  Then $O$ is an
$r$-arc in $\Pi^*$.  Let $S$ denote the set of $q + 1 - r$ points of
$\ell_{\infty}$ not incident with the lines of $\Sigma$.  Extend $O$
to an $(r + 1)$-arc $A$ of $\Pi^*$ by adjoining a point $P$ of $S$.
By Theorem~\ref{t6.1}, $A$ is contained in an oval or hyperoval $B$ of
$\Pi^*$, and $S \subseteq B$ since $A \cup \{ X \}$ is an $(r+2)$-arc
for each point $X$ of $S \setminus \{ P \}$.  As the points of $S$
are collinear, one obtains the contradiction $q + 1 - r \le
2$.

Assume next, by way of contradiction, that $\Pi$  holds an
$r$-net $\Sigma$ with hyperoval $H$ where $r > f(q) - 2$ and either $r
\le q - 3$ or $r \le q - 2$ with $q$ even.  Extend $H$ to an $(r +
3)$-arc $A'$ of $\Pi^*$ by adjoining two points $P$, $Q$ of $S$.  By
Theorem~\ref{t6.1}, $A'$ is contained in an oval or hyperoval $B'$ of
$\Pi^*$.  The upper bound on $r$ guarantees that there is a point $R$
in $B' \setminus A'$.  The $r + 1$ lines of $\Pi^*$ which join $R$ to
points of $H$ all intersect $\ell_{\infty}$ in the $q - 1 - r$ points
of $S \setminus \{ P, Q \}$.  Thus, $r+1 \le q-1-r$; so $(q-2)/2 \ge r
> f(q) - 2$.  Then $q \le 4$, and one obtains the contradiction $r \le
2$.~$\blacksquare$ 

\medskip

Corollary~\ref{t6.2} implies that $q - 2 \notin O^q_d$ for even $q \ge
8$ and  for odd $q \ge 125$ and that $q - 2 \notin H^q_d$ for all even
$q$.  We improve these results in Corollaries~\ref{t6.4} and
\ref{t6.6} below.  

\begin{theorem} \label{t6.3}
\cite[8.6.10]{H}.
Let $q > 3$ be odd.  Then each $q$-arc of $\Pi^* = \Pi^*(GF(q))$ is
contained in a unique oval of $\Pi^*$.
\end{theorem}

\begin{corollary} \label{t6.4}
For odd $q \ge 7 \,$, $q - 2 \notin O^q_d$.
\end{corollary}

\proof Assume, by way of contradiction, that $\Pi := \Pi($GF$(q))$
contains a $(q - 2)$-net $\Sigma$ with oval $O$ for some odd $q \ge
7$.  Adjoin a line $\ell_{\infty}$ to $\Pi$ to form $\Pi^* : =
\Pi^*(GF(q))$.  Let $S := \{P_1$, $P_2$, $P_3\}$ be the set of three
points of 
$\ell_{\infty}$ not incident with lines of $\Sigma$, $T =
\ell_{\infty} \setminus S$, $S_i = S \setminus \{ P_i \}$.  Then $O
\cup S_i$ is a $q$-arc of $\Pi^*$ which, by Theorem~\ref{t6.3}, is
contained in an oval $O_i := O \cup S_i \cup \{ Q_i \}$ of $\Pi^*$.
For given $i$, at least $q - 3$ of the $q - 2$ lines $XQ_i$ with $X
\in O$ contain points of $T$;  i.e., are tangents to $O$ in $\Sigma$.
Since there are only $q - 2$ tangents to $O$ in $\Sigma$, $Q_1$ and
$Q_2$ lie in the intersection of a common pair of tangents, and thus
$Q_1 = Q_2$.  Then $O_1$ and $O_2$ are distinct ovals of $\Pi^*$ whose
intersection is a $q$-arc of $\Pi^*$.  This contradiction of
Theorem~\ref{t6.3} completes the proof of the
corollary.~$\blacksquare$  

\begin{proposition} \label{t6.5}
If an $(n-2)$-net $\Sigma$ is held by an affine plane $\Pi$ of order
$n$, then $\Sigma$ has no hyperovals.
\end{proposition}

\proof Assume, by way of contradiction, that $\Pi$ contains an $(n -
2)$-net $\Sigma$ with hyperoval $H$.  Then $n$ is odd.  Adjoin a
line $\ell_{\infty}$ to $\Pi$ to form a projective plane $\Pi^*$.
Let $S := \{P_1$, $P_2$, $P_3\}$ be the set of three points of 
$\ell_{\infty}$ not incident with lines of $\Sigma$.  Then $H
\cup \{P_1, P_2\} =: O$ is an oval of $\Pi^*$.  By a result of Qvist
(\cite[p.\,148]{DM}), no point of $\Pi^*$ lies on more than two
tangents to $O$.  Since $P_3$ lies on $n - 1$ tangents to $O$, $n \le
3$;  but since $\Sigma$ is an $(n-2)$-net, $n \ge 5$, a contradiction
which yields the asserted conclusion.~$\blacksquare$

\begin{corollary} \label{t6.6}

(i) For every prime power $q$, the integer $q - 2$ is not in $H^q_d$. 

(ii) For every $n \ge 25$, the integer $n - 2$ is not in $H^n$.
\end{corollary}  

\proof By definition, a Desarguesian net of order $q$ is held by an
affine plane of order $q$;  by a well-known theorem of Bruck (see,
e.g., \cite[p.\,714]{BJL}), every $(n-2)$-net of order $n$ with $n \ge
24$ is held by an affine plane of order $n$.~$\blacksquare$

\begin{proposition} \label{t3.12}

(i) $\{ 3 \} \cup H^n \subseteq O^n \subseteq \{ 3, 4,..., n+1 \}$
for all $n$; 

(ii) $4 \in O^n$ if $n \ge 3$ is the order of a projective plane;

(iii) 3, 4, $q \in O^n_d$ if $q$ is a prime power and $n$ is a power
of $q$; 

(iv) $r \in O^q_d$ if $q$ is a prime power and $r$ divides $q+1$ or
$q-1$. 
\end{proposition}

\proof Clearly, $3 \in O^n$ for all $n$, and $3 \in O^n_d$ for prime
powers $n$.  It is well-known that $3 \le r \le n+1$ for every $r$-net
of order $n$ (see, for instance, \cite[Exercise I.7.5]{BJL}); so (i)
follows from Fact~\ref{t2.1}.  If $\Pi$ is an affine
plane of order $n$, let $\ell_1,\ell_2$ be distinct lines in a
parallel class of $\Pi$, let $\ell_3,\ell_4$ be distinct lines in
another parallel class of $\Pi$, and let
$O=(\ell_1\cup\ell_2)\cap(\ell_3\cup\ell_4)$.  Then $O$ is an oval in
a 4-net held by $\Pi$;  so (ii) is valid, and $4 \in O^n_d$ for prime
powers $n$.  The truth of (iv) and the remaining assertions of (iii)
follows immediately from Propositions~\ref{t2.10} and
\ref{t2.5}\,(i).~$\blacksquare$  

\begin{proposition} \label{t6.8} (Drake, Myrvold \cite{DMM})

(i) Every 6-net of order 8 is Desarguesian.

(ii) $5 \in O^8 \cap H^9$;  $6 \in O^9$.
\end{proposition}

\begin{corollary} \label{t3.13}
(i) $O^n_d = O^n = \{r \, \vert \, 3 \le r \le n+1 \}$ for $2 \le n
\le 5$;

(ii) $O^6 = \{3\}$; $O^7_d = O^7 = \{3, 4, 6, 7, 8\}$; 

(iii) $O^8_d = \{3, 4, 7, 8, 9\}$, and $O^8 = O^8_d \cup \{5\}$;

(iv) $O^9_d = \{3, 4, 5, 8, 9, 10\}$, and $O^9 = O^9_d \cup \{6\} \cup
T$ where $T \subseteq \{7\}$;

\end{corollary}

\proof Proposition~\ref{t3.12} (with $n = q$ in part (iii)) determines
the sets $O^n$ and $O^n_d$ with $n \le 5$ and gives $3 \in O^6$, $\{3,
4, 6, 7, 8\} \subseteq O^7_d$, $\{3, 4, 7, 8, 9\} \subseteq O^8_d$,
$\{3, 4, 8, 9, 10\} \subseteq O^9_d$.  Tarry \cite[X.13.1]{BJL} proved
that all $r$-nets of order 6 satisfy $r = 3$;  so $O^6 = \{3\}$.
By Corollary~\ref{t3.7}\,(ii), the integer 5 is in $O^9_d$ but not in
$O^7_d \cup O^8_d$.  By Corollary~\ref{t6.2}\,(i), we have $6
\not\in O^8_d$. Drake used the 
Norton-Sade list of Latin squares of order 7 to observe \cite[Table
1]{D} that every $r$-net of order 7 with $r \ge 5$ can be extended to
an affine plane of order 7.  R. C. Bose and K. R. Nair as 
well as M. Hall \cite[p.\,169]{DK} have proved that the only plane of
order 7 is the Desarguesian one.  Thus, all 5-nets of order 7 are
Desarguesian, so $5 \not\in O^7$.  By Proposition~\ref{t6.8}, $6
\not\in O^8$; but $5 \in O^8$, and $6 \in O^9$.  By Theorem~\ref{t4.6}
and Corollary~\ref{t6.4}, neither 6 nor 7 is in
$O^9_d$.~$\blacksquare$ 

\begin{corollary} \label{t6.10}
(i) $H^2_d = \{3\}$; $H^4_d = \{3, 5\}$; $H^8_d = \{3, 7, 9\}$;
$H^{13}_d \subseteq \{9\}$; $H^{16}_d = \{3, 5, 7, 15, 17\} \cup S$
where $S \subseteq \{ 9 \}$;

(ii) $H^n_d$ is the empty set if $2, 4, 8, 13 \ne n < 16$;

(iii) $H^n = H^n_d \cup \{3\}$ for $8 \ge n \ne 3$; $H^3$ is empty.

(iv)  $\{3, 5\} \subseteq H^9 \subseteq \{3, 5, 7\}$.
\end{corollary}

\proof The sets $H^n$ with $n \le 8$ are given in \cite[Proposition
4.5]{D1}, so it suffices to prove (i), (ii) and (iv).  If $n$ is not a
prime 
power, there is no Desarguesian plane of order $n$; so $H^n_d$ is the
empty set.  By Proposition~\ref{t3.12}\,(i), each $H^n_d$ consists of
odd integers $r$ with $3 \le r \le n + 1$.  By Proposition~\ref{t3.3},
$3 \in H^n_d$ if and only if $n$ is a power of 2.  For $n \le 16$,
Theorem~\ref{t8.3} implies that $5 \in H^n_d$ if and only if $n = 4$
or 16; and Theorem~\ref{t5.1} implies that $7 \in H^n_d$ if and only
if $n = 8$ or 16.  By Proposition~\ref{t2.5}\,(ii) and
Proposition~\ref{t2.4}, $9 \in H^8_d$; and $15, 17 \in H^{16}_d$.  For
every $n$, \cite[Fact 2.3]{D1} yields $n \not\in H^n$;  in particular,
$9 \not\in H^9_d$, $11 \not\in H^{11}_d$, and $13 \not\in H^{13}_d$.
By Corollary~\ref{t6.6}, $9 \not\in H^{11}_d$ and $11 \not\in
H^{13}_d$; by Corollary~\ref{t6.2}\,(ii), 11, $13 \not\in
H^{16}_d$.  By \cite[Corollary 3.1]{D1}, $3 \in H^9$; and $5 \in H^9$
by Proposition~\ref{t6.8}.~$\blacksquare$


\begin{thebibliography}{99}

\bibitem{BJL}  Th.\ Beth, D. Jungnickel and H. Lenz,
{\em Design Theory}, second edition, in Encyclopedia of Math.\ and its
Applications, vol.\ 69, Cambridge Univ.\ Press, Cambridge, (1999).

\bibitem{CA} P. J. Cameron, {\em Parallelisms of Complete Designs}, in
London Math.\ Soc.\ Lecture Note Series, vol.\ 23, Cambridge
Univ.\ Press, Cambridge, (1976).

\bibitem{C} W. Cherowitzo, Hyperovals in Desarguesian planes: an
update, {\it Discrete Math.}, Vol.\ 155 (1996) pp.~31--38.

\bibitem{DM} P. Dembowski, {\em Finite Geometries}, Ergebnisse der
Math. u. ihrer Grenz-gebiete, vol.\ 44, Springer-Verlag, New York,
(1968). 

\bibitem{DK} J. D\'enes and A. D. Keedwell, {\em Latin Squares and
their Applications}, The English Universities Press, London, England,
(1974). 

\bibitem{D} D. A. Drake, Maximal sets of Latin squares and partial
transversals, {\it J. Statistical Planning and Inference}, Vol.\ 1
(1977) pp.~143--149.

\bibitem{D1} D. A. Drake, Hyperovals in nets of small degree, {\it
J. Combinatorial Designs}, to appear.

\bibitem{D2} D. A. Drake, A characterization of the Desarguesian
affine planes which hold 7-nets with ovals, preprint, Univ.\ Florida,
(2002) pp.~1--6.

\bibitem{DMM} D. A. Drake and W. Myrvold, The non-existence of maximal
sets of four mutually orthogonal Latin squares of order 8, preprint
(2002). 

\bibitem{E} R. Evans, H. Hollmann, Ch.\ Krattenthaler and Q. Xiang,
Gauss sums, Jacobi sums, and $p$-ranks of cyclic difference sets,
{\it J. Combinatorial Theory, Series A}, Vol.\ 87 (1999) pp.~74--119. 

\bibitem{H} J. W. P. Hirschfeld, {\em Projective Geometries
over Finite Fields}, Oxford University Press, Oxford, (1979).

\bibitem{L} T. Y. Lam, {\em A First Course in Noncommutative Rings},
second edition, Graduate Texts in Mathematics, Springer-Verlag, New
York, (2001).

\bibitem{S} J.-P. Serre, {\em A Course in Arithmetic},
Springer-Verlag, New York, (1973).

\bibitem{Simm} G. J. Simmons, Sharply focused sets of lines on a conic
in $PG(2,q)$, {\it Congressus Numerantium}, Vol. 73 (1990)
pp.~181--204. 

\end{thebibliography}
\end{document}